\documentclass[12pt]{amsart}

\def\ca{{\mathcal A}}

\def\Mp{M_{p\times p}}
\def\Np{N_{(p-1)\times (p-1)}}

\newcommand{\comment}[1]{}
\input{amssym.def}
\input{amssym.tex}
\newtheorem{theorem}{Theorem}[section]

\newtheorem{corollary}[theorem]{Corollary}
\newtheorem{remark}[theorem]{Remark}
\newtheorem{definition}[theorem]{Definition}

\title[Canonical Representatives for Conjugacy Classes]{Canonical
Symplectic Representations for Prime Order Conjugacy Classes of the
Mapping-class Group}
\thanks{Partially supported by grants from the NSA and the Rutgers Research Council and by Yale University}
\author{Jane Gilman}
\address{Mathematics Department, Rutgers University, Newark, NJ 07102}%
\email{gilman@rutgers.edu}%

\begin{document}

\begin{abstract}

In this paper we find a unique normal form for the symplectic matrix
representation of the conjugacy class of a prime order element of
the mapping-class group. We find a set of generators for the
fundamental group of a surface with a conformal automorphism of
prime order which reflects the action the automorphism in an optimal
way. This is called an {\sl adapted} homotopy basis and there is a
corresponding {\sl adapted presentation}. We also give a necessary
and sufficient condition for a prime order symplectic matrix to be
the image of a prime order element in the mapping-class group.
\end{abstract}
\maketitle
\section{Introduction}

This is the second of two papers on prime order conformal
automorphisms of compact Riemann surfaces of genus $g \ge 2$. The
mapping-class group of a surface $S$ of genus $g$, denoted by
$MCG(S)$,  is the set of homotopy classes of self-homeomorphisms of
$S$. It is well known that the action of a homeomorphism $h$ of a
surface of genus $g$ on a {\sl canonical homology basis}, a homology
basis with certain intersection properties, will give a symplectic
matrix, $M_{h_{CAN}}$. That is, $M_{h_{CAN}}$ is an element of
$Sp(2g,\mathbb{Z})$, the set of $2g \times 2g$ integer valued
matrices that preserve the symplectic form, the  non-degenerate
skew-symmetric
bilinear quadratic form with matrix, $J = \left(%
\begin{array}{cc}
  0 & I_g \\
  -I_g & 0 \\
\end{array}%
\right)$ where $I_g$ denotes the $g \times g$ identity matrix.  The
intersection numbers for the curves in a canonical basis is
precisely given by the matrix $J$. We fix a canonical homology basis
and let $\pi$ denote the surjective map $\pi: MCG(S) \rightarrow
Sp(2g,\mathbb{Z})$.

In \cite{G4, JGind} we found that corresponding to any conformal
automorphism $h$  of prime order, there exists a special set of
generators for the first homology of the surface that reflected the
action of the homeomorphism $h$ in an optimal way. Such a basis was
termed an {\sl adapted basis} or equivalently  a {\sl basis adapted
to $h$}. In particular the matrix action of the automorphism on the
basis had a very simple form. The curves in this homology basis did
not form a canonical homology basis. The intersection matrix for the
adapted basis was found in \cite{GPint}.

In particular, we have ordered pairs of matrix representations and
intersection matrices $(M_{h_{CAN}}, J)$ and $(M_{\mathcal{A}},
J_{\mathcal{A}})$  where $M_{\mathcal{A}}$ is the matrix of the
action on an adapted basis and $J_{\mathcal{A}}$ the intersection
matrix for that basis. For each $0 < t \le 2g +2$ and for each prime
$p \ge 2$, there is a set of $t$ integers $(n_1,...,n_t)$ with $0 <
t \le p-1$ with $\Sigma_{j=1}^tn_j \equiv 0 (p)$ that determines the
conjugacy class of $h$. Here $\equiv ... (p) $ denotes equivalence
modulo $p$. These integers are reflected in each ordered pair.
However, in the second pair, these integers are reflected in the
intersection matrix and in the first pair, reflected in the matrix
of the action of the homeomorphism.

In this paper we present an algorithm to obtain $M_{h_{CAN}}$ from
$M_{\mathcal{A}}$. The input to the algorithm is the $t$-tuple of
integers. The output is a canonical representative for the conjugacy
class of $h$ in $Sp(2g, \mathbb{Z})$. Thus the main result of this
paper is to describe a canonical representative for each prime order
conjugacy class in $Sp(2g, \mathbb{Z})$ that is the representation
of a prime order element of the mapping-class group. All other prime
order symplectic matrices are representations of infinite order
mapping-classes.

The organization of this paper is as follows: In Part 1, we fix
notation and summarize terms and prior results; in Part 2 we develop
the tight surface symbol reduction algorithm, and in  Part 3 we
review the basics of the Schreier-Reidemeister rewriting process. In
Part  4, using new calculations we obtain an {\sl adapted
presentation of the fundamental group} corresponding to a conformal
or finite order mapping class. The main theorems (theorem
\ref{theorem:BIG} and theorems \ref{theorem:algbig},
\ref{theorem:algbignext} and \ref{theorem:full}) are proved there in
Part 4 and a detailed example is worked out by hand (section
\ref{section:ex}).

There is renewed interest in mapping-class groups in regard to
computing on Riemann surfaces. In earlier work we noted that in the
abstract there was in theory a method  for converting the
 adapted basis to a canonical homology basis and the adapted matrix
 representation to a symplectic matrix representation, but the formula involved
 was unwieldy.  However, we have recently observed that while the
 formula
is messy for hand computations, it can be implemented as an
algorithm whose complexity is bounded by $O(g^{10})$.

 \tableofcontents

\part{Summary of terms and prior results }
\section{Preliminaries} \label{section:preliminaries}
\subsection{Notation and Terminology}
 We let  $h$ be a conformal automorphism of
a compact Riemann surface $S$ of genus $g \ge 2$. Then $h$ will have
a finite number, $t$,  of fixed points. We let $S_0$ be the quotient
of $S$ under the action of the cyclic group generated by $h$ so that
$S_0 = S/\langle h \rangle$ and let $g_0$ be its genus. If $h$ is of
prime order $p$ with $p \ge 2$ , then the Riemann-Hurwitz relation
shows that
 $2g = 2pg_0 + (p-1)(t-2)$. If $p=2$, of course, this implies that $t$ will be
 even.

\subsection{Equivalent Languages} \label{section:tor}

 We emphasize that $h$ can be thought of in
a number or equivalent ways using different terminology. For a
compact Riemann surface of genus $g \ge 2$, homotopy classes of
homeomorphisms of surfaces are the same as isotopy classes.
Therefore, $h$ can be thought of as a representative of a homotopy
class or an isotopy class. Further, every isotopy class of finite
order contains an element of finite order so that $h$ can be thought
of as a homeomorphism of finite order. For every finite order
homeomorphism of a surface there is a Riemann surface on which its
action is conformal. A conformal homeomorphism of finite order up to
homotopy is finite. We use the language of conformal maps, but
observe that all of our results can be formulated using these other
classes of homeomorphisms.

We remind the reader that the mapping-class group of a compact
surface of genus $g$ is also known as the Teichm{\"u}ller Modular
group or the Modular group, for short. We write $MCG(S)$ or
$MCG(S_g)$ for the mapping class group of the surface $S$ using the
$g$ when we need to emphasize that $S$ is a compact surface of genus
$g$. The Torelli Modular group or the Torelli group for short,
${\mathcal T}(S)$ is homeomorphisms of $S$ modulo those that induce
the identity on homology and the homology of a surface is the
abelianized homotopy. There is surjective map $\pi$ from the
mapping-class group onto $Sp(2g,\mathbb{Z})$ that assigns to a
homeomorphism the matrix of its action on a canonical homology basis
(see \ref{section:homology}). The induced map from the Torelli group
to the symplectic group is an isomorphism.

It is well known that the map $\pi$ when restricted to elements of
finite order is an isomorphism. An even stronger result can be found
in \cite{G4, JGind}.

Since  $h$ can be thought of as a finite representative of a finite
order mapping-class, we will always treat it as finite.  For ease of
exposition we use the language of a conformal maps and do not
distinguish between a homeomorphism that is of finite order or that
is of finite order up to homotopy or isotopy, a finite order
representative for the homotopy class, a conformal representative
for the class or the class itself. That is, we do not use different
notation to distinguish between the topological map, its homotopy
class or a finite order representative or a conformal
representative.

For ease of exposition in what follows we first assume that  $t
> 0$. We treat the case $t=0$ separately in section
\ref{section:t=0}.

\subsection{Conjugacy Invariants for prime order mapping classes or
conformal automorphisms} \label{section:mandn}

 Nielsen showed that the conjugacy class
of  $h$ in the mapping-class group is determined by the set of $t$
non-zero integers $\{n_1, ...., n_t\}$ with $0 < n_i < p$ where
$\Sigma_{i=1}^t n_i \equiv 0 \; (p)$.

 Let $m_j$ be the number of $n_i$
equal to $j$. Then we have $\Sigma_{i=1}^{p-1} i \cdot m_i \equiv 0
\; (p)$ (see \cite{G3} for details) and the conjugacy class is also
determined by the $(p-1)$-tuple, $(m_1,...,m_{p-1})$.

Topologically we can think of $h$ as a counterclockwise rotation by
an angle of  ${\frac{2\pi \cdot s_i}{p}}$ about the fixed point
$p_i, i = 1, ..., t$ of $h$. We call the $s_i$ the {\sl rotation
numbers}. The $n_i$ are the {\sl complementary rotation numbers},
that is, $0 < s_i < p$ with $s_in_i \equiv 1 \; (p)$.

\subsection{Homology}\label{section:homology}
We recall the following facts about Riemann surfaces.

The homology group of a compact Riemann surface of genus $g$ is the
abelianized homotopy. Therefore, a homology basis for $S$ will
contain $2g$ homologously independent curves. Every surface has a
{\sl canonical homology basis}, a set of $2g$ simple closed curves,
$a_1,...,a_g; b_1,...,b_g$ with the property that for all $i$ and
$j$,  $a_i \times a_j = 0$,
 $b_i \times b_j = 0$ and $a_i \times
 b_j = \delta_{ij} = -b_j \times  a_i$ where $\times$ is the algebraic intersection number
 and $\delta_{ij}$ is the Kronecker delta.

\subsection{Convention for actions on curves, homotopy, and homology}
A homeomorphism $h$ of a surface induces an action on the
fundamental group of the surface (mapping the fundamental group with
base point $x_0$ to the one with bases point $h(x_0)$) and thus acts
as an (outer) automorphism of the  fundamental group of $S$. It also
induces an automorphism of the first homology group of $S$. If
$\gamma$ is any curve on $S$ or any representative of a free
homotopy class or homology class on $S$, we adopt the convention
that $h(\gamma)$ denotes either the image curve or the free homotopy
class or homology class of the image. It will be clear from the
context which we mean.

\section{Prior results and definitions: matrices and homology bases} \label{section:adapted}

 Roughly speaking a homology basis for $S$ is {\it adapted to $h$} if
it reflects the action of $h$ in a simple manner:  for each curve
$\gamma$ in the basis either all of the images of $\gamma$ under
powers of $h$ are also in the basis or the basis contains all but
one of the images of $\gamma$ under powers of $h$ and the omitted
curve is homologous to the negative of the sum of the images of
$\gamma$ under the other powers of $h$.

To be more precise

\begin{definition} \label{definition:adapt}  A homology basis
for $S$ is adapted to $h$ if for each $\gamma_0$ in the basis there
is a curve $\gamma$ with $\gamma_0 =h^k(\gamma)$ for  some integer
$k$ and either
\begin{enumerate}
\item  $\gamma, h(\gamma), ...h^{p-1}(\gamma)$ are all in the basis,
or
\item $\gamma, h(\gamma), ...h^{p-2}(\gamma)$ are all in the basis
and

$h^{p-1}(\gamma) \approx^h -(h(\gamma) +  h(\gamma)+ ...
+h^{p-2}(\gamma))$.

Here $\approx^h$ denotes is homologous to.

\end{enumerate}

\end{definition}

\subsection{Existence of adapted homology bases}
\label{section:exist}
 It is known that
\begin{theorem} {\label{theorem:ad}}
 {\rm (\cite{G4}, \cite{JGind})} There is a homology basis adapted to $h$. In
particular, if  $g\ge 2$, $t \ge 2$,  $g_0$ are as above, then the
adapted basis has $2p \times  g_0$ elements of type (1) above and
$(p-1)(t-2)$ elements of type (2).
\end{theorem}
and thus it follows that
\begin{corollary} {\rm (\cite{G4})} \label{corollary:corad}
Let  $M_{\ca}(h)$ denote the {\sl adapted matrix of $h$}, the matrix
of the action of $h$  with respect to an adapted basis. Then
$M_{\ca}(h)$ will be composed of diagonal blocks, $2g_0$ of which
are $p \times p$ permutation matrices with $1$'s along the super
diagonal and $1$ in the leftmost  entry of the last row and $t$ are
$(p-1) \times ( p-1)$ matrices with $1$'s along the super diagonal
and all entries in the last row $-1$.
\end{corollary}

\begin{remark} We adopt the following convention.
When we pass from homotopy to homology, we use the same notation for
the homology class of the curve as for the curve or its homotopy
class, but write $\approx^h$ instead of $=$. It will be clear from
the context which we mean.
\end{remark}

\subsection{Intersection Matrix for an Adapted Homology Basis}
\label{section:intersection} So far  information about $M_{\ca}$
seems to depend only on $t$ and not upon the $(p-1)$-tuple
$(m_1,...,m_{p-1})$ or equivalently, upon the set of integers
$\{n_1, ..., n_t\})$ which determines
 the conjugacy class of $h$ in the mapping-class group.  However,
while the $2pg_0$ curves can be extended to a canonical homology
basis for $h$, the rest of the basis can  not and its intersection
matrix, $J_{\ca}$ depends upon these integers.

In \cite{GPint} the intersection matrix for the adapted basis was
computed.

The adapted basis consisted of the curves of type (1):
$$\{A_w, B_w,  w=1,..., g_0 \} \cup \{ h^j(A_w), h^j(B_w),
j=1,...,p-1 \}$$ and (some of) the curves of type (2):
$$X_{i,v_i}, h^j(X_{i,v_i}),\;  i=1,...,(p-1),\; j =1, ..., p-2,\;
v_i=1,...,u_i.$$

A lexicographical order is placed on $X_{i,v_i}$ so that $(r,v_r) <
(s, v_s)$ if and only if $r < s$ or $r=s$ and $v_r < v_s$. The $t-2$
curves $X_{s,v_s}$ with the largest subscript pairs are to be
included in the homology basis. Let $\hat{s}$ be the smallest
integer $s$ such that $u_s \ne 0$ and let $\hat{q}$ be chosen so
that $\hat{q} \hat{s} \equiv 1 (p)$. For any integer $v$ let $[v]$
denote the least non-negative residue of $\hat{q} v$ modulo $p$.
Thus the integer $[v]$ satisfies $0 \le [v] \le p-1$ and $\hat{s}
\times [v] \equiv v $ mod(p).

\begin{theorem} {\label{theorem:GPint}}
{\rm (\cite{GPint})}  If $(u_1,...,u_{p-1})$ determines the
conjugacy class of $h$ in the mapping-class group, then the surface
$S$ has a homology basis consisting of:

\begin{enumerate}
\item  $h^j(A_w), h^j(B_w)$  where $1 \le w \le g_0,  0 \le j \le p-1$.

\item $h^k(X_{s,{v_s}})$   where $0 \le k \le p-2$ and for all pairs $(s,v_s)$
with $1 \le s \le p-1, 1 \le v_s \le u_s$ except that the two
smallest pairs are omitted.

 The intersection numbers for the elements of the
adapted basis are given by
\begin{enumerate}
\item $h^j(A_w) \times h^j(B_w) = 1$

\item If $(r,v_r) < (s,v_s)$, then

$ {h^0(X_{r,{v_r}}) \times h^k(X_{s,{v_s}})} = \left\{
\begin{array}{ll}
        1 & {\mbox{ if }} \;\; [k] < [r] \le [k+s]    \\
       -1 & \mbox{ if } \;\; [k+s] < [r] \le [k]
                                    \end{array}
                                         \right. $

$ {h^0(X_{s,{v_s}}) \times h^k(X_{s,{v_s}})}
                  =  \left\{  \begin{array}{ll}
                        1 & {\mbox{ if } \;\; [k] \le [s] <  [k+s]} \\
\                       -1 & {\mbox{ if } \;\; [k+s] < [s]  < [k]}
                   \end{array}
                      \right.  $
\end{enumerate}

\item All other intersection numbers are $0$ except for those that
following from the above by applying the identities below to
arbitrary homology classes $C$ and $D$. \subitem $ C \times D = - D
\times C$

\subitem $ h^j(C) \times h^k(D)  = h^0(C) \times h^{k-j}(D)$, ($k-j$
reduced modulo $p$.)
\end{enumerate}
\end{theorem}

\subsection{Matrix forms \label{section:matrixforms}}

 We can write the results of theorems
\ref{theorem:ad} and \ref{theorem:GPint} and corollary
\ref{corollary:corad} in an explicit matrix form. To do so we fix
notation for some matrices. We will use the various explicit forms
in subsequent sections.

We let $M_{{\tilde{\mathcal{A}}}}$ denote the matrix of the action
of $h$ on an adapted basis and $J_{{\tilde{\mathcal{A}}}}$ be the
corresponding intersection matrix. Further, we let $M_{h_{CAN}}$ be
the matrix of the action of $h$ on a canonical homology basis. The
corresponding intersection matrix is denoted by $J_{h_{CAN}}$ or
$J$.

However, we prefer to replace $J$  by the following block matrix
where $q= {\frac{(p-1)(t-2)}{2}}$ and where for any integer positive
integer $d$,  $I_d$ denotes the $d \times d$ identity matrix.

$$J_{h_{CAN}} = \left(
\begin{array}{cccc}
0 & I_{pg_0} & 0 & 0\\
-I_{pg_0} & 0 & 0 & 0\\
0 & 0 & 0 & I_{q} \\
0 & 0 & -I_q & 0 \\
\end{array} \right). $$

 We denote the
$p\times p$ permutation matrix by
$$M_{p\times p}  =  \left(
\begin{array}{ccccccc}
0 & 1 & 0 & 0  &\ldots &0 & 0\\
0 & 0 & 1 & 0  &\ldots &0 & 0\\
0 & 0 & 0 & 1  &\ldots &0&  0\\
& & & & \ddots \\
0 & 0 & 0 & 0  &\ldots &1& 0\\0 & 0 & 0 & 0  &\ldots& 0 & 1\\
1 & 0 & 0 & 0  &\ldots &0 & 0\\
\end{array} \right),$$ the $(p-1)\times (p-1)$ non-permutation matrix of the theorem by
$$N_{(p-1)\times(p-1)}  = \left(
\begin{array}{rrrrrrr} 
0 & 1 & 0 & 0  &\ldots &0 & 0\\
0 & 0 & 1 & 0  &\ldots &0 & 0\\
0 & 0 & 0 & 1  &\ldots &0&  0\\
& & & & \ddots \\
0 & 0 & 0 & 0  &\ldots &1& 0\\0 & 0 & 0 & 0 & \ldots &0 & 1\\
-1 & -1 & -1 & -1&  \ldots& -1 & -1\\
\end{array} \right) $$

Thus we have the  $2g_0p \times 2g_0p$ block matrix

$$M_{{\mathcal{A}}_{2g_0,  p\times p}} = \left(
\begin{array}{ccccccc}
\Mp & 0& 0  &\ldots &0 & 0\\
0 & \Mp & 0   &\ldots &0 & 0\\
& & & \ddots \\
0 & 0 &  0  &\ldots &\Mp& 0
\\0 & 0  & 0  &\ldots& 0 & \Mp
\\
\end{array} \right) $$

and the $(t-2)\cdot (p-1) \times (t-2) \cdot (p-1)$ block matrix
$N_{{\mathcal{A}}_{(t-2), (p-1)\times (p-1)}}$

 $$  \left(
\begin{array}{ccccccc}
\Np & 0& 0   &\ldots &0 & 0\\
0 & \Np &  0  &\ldots &0 & 0\\
& & & \ddots& \\
0 & 0 & 0   &\ldots &\Np& 0
\\0 & 0 &  0  &\ldots& 0 & \Np\\
\end{array} \right) $$ so that the $2g \times 2g $ matrix
$M_{\mathcal{A}}$ breaks into
blocks and can be written as
$$M_{\mathcal{A}} = \left(
\begin{array}{cc}
M_{{\mathcal{A}}_{2g_0, p\times p}} &  0 \\
 0&  N_{{\mathcal{A}}_{(t-2),  (p-1)\times (p-1)}}
\\
\end{array} \right) $$ where the blocks are of appropriate size.
The basis can be rearranged so that  $2g \times 2g $ matrix
$M_{\mathcal{\tilde{A}}}$ corresponding to the rearranged basis
breaks into blocks
$$M_{\mathcal{\tilde{A}}} = \left(
\begin{array}{ccc}
M_{{\mathcal{A}}_{g_0, p\times p}} &  0 & 0\\
0 & M_{{\mathcal{A}}_{g_0, p\times p}} &  0 \\
0 & 0&  N_{{\mathcal{A}}_{(t-2),  (p-1)\times (p-1)}}
\\
\end{array} \right) $$ Here
 the submatrix $$\left( \begin{array}{cc}
M_{{\mathcal{A}}_{g_0, p\times p}} &  0 \\
0 & M_{{\mathcal{A}}_{g_0, p\times p}}  \\
\end{array} \right) $$ is a symplectic matrix. We obtain the corollary

\begin{corollary} Let $S$ be a compact Riemann surface of genus $g$ and
assume that $S$ has a conformal automorphism $h$ of prime order $p
\ge 2 $. Assume that $h$ has $t$ fixed points where $t \ge 2$. Let
$S_0$ be the quotient surface $S_0 = S / \langle h \rangle $ where
$\langle h \rangle$ denotes the cyclic group generated by $h$ and
let $g_0$ be the genus of $S_0$ so that  $2g = 2pg_0 + (t-2)(p-1)$.
\vskip .1in
 There is a homology bases on which the action of $h$
is given by the $2g \times 2g$ matrix $M_{\mathcal{\tilde{A}}}$.
\vskip .1in
 The matrix $M_{\mathcal{\tilde{A}}}$ contains a $2g_0p
\times 2g_0p$ symplectic submatrix, but is not a symplectic matrix
except in the special case $t = 2$.
\end{corollary}

\begin{remark} We note that if $p=2$,
$\Mp$ reduces to $ \left(
\begin{array}{cc}
0& 1\\
1 &  0
\\
\end{array} \right) $ and $\Np$ to the $1\times1$ matrix $-1$.
\end{remark}

The point here is that while two automorphisms with the same number
of fixed points will have the {\it same}  matrix representation with
respect to an adapted basis, the intersection matrices will not be
the same and, therefore, the corresponding two matrix
representations in the symplectic group will not be conjugate.

We seek an algorithm to  replace $M_{{\tilde{\mathcal{A}}}}$ by a
the symplectic matrix $M_{h_{CAN}}$ by replacing the $2q \times 2q$
submatrix $N_{(t-2),(p-1)\times(p-1)}$ by a symplectic matrix of the
same size. We will call this matrix $N_{symp{\tilde{\mathcal{A}}}}$.

We note that $J_{{\tilde{\mathcal{A}}}}$ is of the form
$$\left(
\begin{array}{cccc}
 0& I_{pg_0}& 0   & 0\\
-I_{pg_0} & 0 & 0 & 0  \\
0 &0 &B_1 &B_2\\
0& 0& B_3& B_4 \\
\end{array} \right) $$
where the blocks $B_i$ are of the appropriate dimension and we let
$B$ denote the $2q \times 2q$ matrix
$$\left(
\begin{array}{cc}
B_1 &B_2\\
 B_3& B_4 \\
\end{array} \right). $$

We also note that $M_{h_{CAN}}$ will break up into

$$M_{h_{CAN}} =
\left(
\begin{array}{ccc}
M_{{\mathcal{A}}_{g_0, p\times p}} &  0 & 0\\
0 & M_{{\mathcal{A}}_{g_0, p\times p}} &  0 \\
0 & 0&  N_{h_{CAN}}\\
\end{array} \right) $$ where $N_{h_{CAN}}$ is a $(t-2)(p-1)\times (t-2)(p-1)$ matrix.

We emphasize that our goal is to find a canonical form for
$N_{symp{\tilde{\mathcal{A}}}}= N_{h_{CAN}}$ and an algorithm that
produces it from the conjugacy class data. This is found in section
\ref{section:LOOK}. We illustrate this with a detailed example in
section \ref{section:ex}.

\part{The Surface Symbol Algorithm}

\section{The Surface symbol reduction algorithm}

 For a compact Riemann
surface of genus $g \ge 2$, a presentation for the fundamental group
with one defining
 relation determines a {\sl surface symbol}. There is a standard way to convert any surface
 symbol to a normal form \cite{Bers, MacBDun} where it is written as the product of the minimal number
  of commutators.  Here we follow the details and notation from
  Springer (section 5.5 of  \cite{SPR}).

 If we have a compact Riemann surface, we obtain the {\sl surface
 symbol} homeomorphic to the surface under identifications by
cutting along the generators of the fundamental group and labeling
the sides appropriately. Let ${\mathcal{P}}$ denote the polygon
obtained by cutting along the curves.
 Since ${\mathcal{P}}$ is
 obtained by cutting along curves in the fundamental group, each
 {\sl side} occurs once in the positive direction and once in the
 negative direction as one moves along the boundary in the positive (counterclockwise) direction.
  The surface symbol is the sequence of {\sl
 sides} that occur. That is it is a sequence of the form $abc \ldots
 a^{-1}\ldots c^{-1}\ldots rst b^{-1}\ldots$.

We define
\begin{definition} \label{definition:properties}
\begin{enumerate}
\item
A polygon is {\sl evenly worded} if for each edge   $a$ that occurs
$a^{-1}$ also occurs.
\item A pair of edges $a$ and $b$ are {\sl linked}  if they appear in
the symbol in the  order $ \ldots a \ldots b \ldots a^{-1}
\ldots b^{-1} \ldots $.
\item
The polygon is {\sl fully linked} if each edge occurring in the
relation is linked to another unique distinct  edge.
\item The edges  $a$ and $b$ occurring in the polygon are {\sl
optimally linked} if  the symbol is of the form
$W_0aW_1bW_2a^{-1}W_3b^{-1}W_4$ where the $W_i, i=0,...,4$ are words
in the edges not involving $a^{\pm1} $ or $b^{\pm 1}$.

\item  A pair of edges $m$ and $n$ are {\sl tightly linked} if they appear
in the symbol in the  order $mnm^{-1}n^{-1}$.
\item A polygon is {\sl tightly worded} if each edge or its inverse is tightly linked to a
another edge in the symbol.
\end{enumerate}
\end{definition}

Note that a symbol is tightly worded when it can be written as a
product of commutators and the product of $g$ commutators is the
normal form for a compact surface of genus $g$. Thus it corresponds
to a canonical homology basis. Also note that for the symbol we will
be tightly linking all linked edges are optimally linked because
each edge and its inverse occur only once in the symbol and the
steps we take will not change that. Finally, a surface symbol will
always have an even number of edges since each edge and its inverse
will occur.

\begin{remark}
More generally, one can consider arbitrary polygons with identified
sides labeled as $x$ and $x^{-1}$. In that case the polygon may have
more than one equivalence class of identified vertices. There is an
algorithm for replacing the surface symbol by another surface symbol
in a manner that increases the number of vertices in one equivalence
class and decreases the number of vertices in the other until one
obtains a symbol with only one equivalence class of vertices
\cite{SPR}. Since we are taking our surface symbol from the
fundamental group with fixed base point, our surface symbols only
have one equivalence class of vertices.
\end{remark}

In \cite{SPR} it is shown that if the symbol has one vertex up to
identification, which it does if one takes the fundamental group
with fixed base point, then each edge is linked with some other
edge.

A surface symbol may posses one or more of the properties described
in definition \ref{definition:properties}.

Our goals is to begin with a surface symbol and obtain a tightly
worded polygon.

\subsection{How to reduce a linked
polygon}\label{section:reducelink} Let $a$ and $b$ be  an optimally
linked pair of edges with linking words $W_0, W_1,W_2, W_3, W_4$ so
that the symbol is
 $W_0 a W_1 b W_2a^{-1} W_3 b^{-1} W_4...$. Then one can replace $a$ and $b$ by new words $M$ and $N$
 and obtain a new polygon where the new words $M$ and $N$  are tightly
 linked by a standard cutting and pasting
 replacements. In terms of the surface symbol the words $M$ and $N$ now represent edges.

 We apply this cutting and pasting procedure in the following informal algorithm.
 This algorithm does more than tightly link the surface symbol. The algorithm input includes
 the adapted matrix and its output includes a symplectic matrix representation.
We let $\lambda$ represent the empty word. \vskip .1in
{\centerline{\bf Ordered recipe for tight linking to a symplectic
matrix.}

\begin{description}
\item [Step 1]{\sl Initialize} Assume that the polygon ${\mathcal{P}}$ has
edges $e_1,..., e_m$ where $m$ is divisible by $4$. Thus
${\mathcal{P}} = e_1e_2 \cdots e_m$. Further assume that polygon
${\mathcal{P}}$ is evenly worded, that it has only one
equivalence class of vertices and that it is fully and optimally
linked. The initial matrix is $M_{\mathcal{A}}$. Input
${\mathcal{P}}$ and $M_{\mathcal{A}}$. Set ${\mathcal{M}} =
M_{\mathcal{A}}$. Set $\mathcal{Q} = \lambda$

\item[Step 2] {\bf DO}
\begin{itemize}
    \item (Free reduction)    If $aa^{-1}$ or $a^{-1}a$ occur in ${\mathcal{P}}$ they can be removed from the
symbol as long as the symbol has at least one other letter. Delete
$a$ and $a^{-1}$ from the set of edges and renumber the remaining
edges. Replace ${\mathcal{P}}$ by ${\mathcal{P}}$ with $aa^{-1}$
deleted. Repeat this step wherever possible, until no repeated
symbols of the form $xx^{-1}$ or $x^{-1}x$ occur. ({\sl Note: we
will see that in the case where we apply this algorithm we will
always have at least one other letter.})

    \item (Linking and truncation)
Begin with left most occurring letter in the surface symbol,
${\mathcal{P}}$. Let $a$ be this letter. Let $b$ be the {\sl first}
edge optimally linked with $a$. Write ${\mathcal{P}} =$
$aW_1bW_2a^{-1}W_3b^{-1}W_4$ where $W_1, W_2, W_3$ and $W_4$ are
chosen accordingly. Define $M= aW_1bW_2a^{-1}$ and $N=
W_3W_2a^{-1}$.
\\({\sl Note: ${\mathcal{P}}$ can now be
written as $[M,N]W_3W_2W_1W_4$.})\\ Preserving cyclical order rename
the edges other than \\$a, a^{-1}, b, \mbox{ and } b^{-1}$.  That
is, define $e_1',...,e_m'$ to be the remaining edges that occur in
$W_3W_2W_1W_4$.

({\sl Note: We have the initial ordered basis $e_1,...,e_m$ and new
ordered basis $M,N,M^{-1},N^{-1}, e_1',...,e_{m'}'$. In terms of the
surface symbol the words $M$ and $N$ now represent edges. Note also
that $a, b, a^{-1},$ and $b^{-1}$ are each one of the original
$e_j's$.})

\item (Change of Basis)

Write down the matrix of the change of basis ${\mathcal{B}}$.
Namely, ${\mathcal{B}}$ is the identity matrix except that the
columns corresponding to $a$ and $a^{-1}$ are replaced by the
columns whose entries are determined  by $M$ and similarly for $b$
and $b^{-1}$ and $N$. Compute
${\mathcal{B}}{\mathcal{M}}{\mathcal{B}}^{-1}$.

\vskip .1in
 Replace the matrix ${\mathcal{M}}$ by its conjugate
${\mathcal{B}}{\mathcal{M}}{\mathcal{B}}^{-1}$.

Set $e_i = e_i', \;\; i=1,...,m'=m-4$.

Set ${\mathcal{Q}} = {\mathcal{Q}} \cdot [{\tilde{a}}, {\tilde{b}}]$

Set ${\mathcal{P}} = W_3W_2W_1W_4$ \vskip .1in

{\centerline{{\bf UNTIL}: ${\mathcal{P}} = \lambda$}}
\end{itemize}
\vskip .1in

\item[Step 3] Output: ${\mathcal{Q}}$ and ${\mathcal{M}} =
M_{h_{CAN}}$. \vskip .05in \noindent  ({\sl Note: the output
consists of the tightly worded surface symbol and the symplectic
matrix.})
\end{description}

\subsection{Analysis of the algorithm}

The input to the algorithm will be the surface symbol
${\hat{\mathcal{Q}}} \cdot {\hat{\hat{\mathcal{LR}}}}$ and the
matrix $M_{\tilde{\tilde{{\mathcal{A}}}}}$ of corollary
\ref{corollary:adsymb}. Alternately the input to the algorithm can
be taken to be (a variant of) the  surface symbol ${\hat{\hat{R}}}$
derived in section \ref{section:LOOK}.
  It has length $m = 4g$. The free reduction step is
$O(m)$. It reduces the number of edges. In our application we will
begin with $4g$ edges that represent a surface of genus $g$. We know
that the final surface symbol must also represent a surface of genus
$g$ so that the final surface symbol will still have $4g$ edges.
Thus we will never encounter a free reduction step (see also remark
\ref{remark:new}). Each repetition of the loop reduces $m$ by $4$ so
the algorithm will stop after at most ${\frac{m}{4}}=g $ steps.

We can actually initialize our algorithm with any ${\mathcal{Q}},
{\mathcal{M}}, \mbox{ and } {\mathcal{P}}$ of the correct sizes as
long as ${\mathcal{Q}}$ is a product of commutators and
${\mathcal{P}}$ evenly worded and fully and optimally linked. The
length of ${\mathcal{P}}$  will be at most $4g$.

In the linking step we need to do a search to find the {\sl first}
$b$ optimally linked to $a$. Since we are searching a string of
length at most $m$ the time here is  $O(m^2)$.

 The matrix of the change of basis is of size $2g \times 2g$ and  will only have four
  columns that differ from the identity matrix, those rows and columns
 that represent replacing  $a$ by the edge labeled  $M$ and  $b$ by the edge now labeled
 $N$.

We are working with free homotopy classes. These are not necessarily
abelian, but when we write down the matrix representation and the
matrix of the change of basis, we can cancel any occurrence of a
generator (or edge) $x$ and $x^{-1}$ in the words $aW_1bW_2a^{-1}$
and $W_3W_2a^{-1}$. The words we are working with are of length at
most $m$. Thus the matrix of the change of basis can be computed in
$O(m^2)$ time. The matrix multiplications and the inversion for the
change of basis step are  $O(m^3)$.

 We conclude
that,  neglecting the size of the integers involved,  the complexity
of the tight linking algorithm in the case of interest is  $O(g^3)$.
The size of the input is taken into account  in section
\ref{section:complex}.

\begin{remark} \label{remark:new} We note that iterating the loop for free reduction and linking is the standard
algorithm for finding the normal form for a surface symbol (section
5.5 of \cite{SPR}). It can be applied to any evenly worded symbol
that comes from the fundamental group of a surface. The symbol does
not need to be of length divisible by $4$. It is guaranteed to give
either a final symbol of the form $aa^{-1}$ or a product of $2n$
commutators for some integer $n$ (Theorem 5-16 of \cite{SPR}). In
our case, since we want to add the truncation and change of basis
steps, we have assumed the divisibility by $4$ to make the
description of these steps simpler. Also we have $2g$ generators for
our fundamental group, those in equation \ref{equation:genssym} of
theorem \ref{theorem:BIG}  or conjugates of those (corollary
\ref{corollary:symbol} or corollary \ref{corollary:adsymb}). Since
each generator and its inverse appear exactly once in our surface
symbol, we have $4g$ edges for the input. We know by construction
that this is a surface symbol for a compact surface of genus $g$.
Therefore, we  must end with $4g$ edges (i.e. $g$ commutators) so in
our implementation of the algorithm, we must never encounter a free
reduction step.
\end{remark}

\part{Review of rewriting basics, prior results and notation}
\section{Schreier-Reidemeister Rewriting and its corollaries} \label{section:SCHR}

If we begin with an arbitrary finitely presented group $G_0$ and a
subgroup $G$, the Schreier-Reidemeister rewriting process tells one
how to obtain a presentation for $G_0$ from the presentation for
$G$. In our case the larger group $G_0$ will correspond to the group
uniformizing $S_0$ and the subgroup $G$ corresponds to the group
uniformizing $S$.

In particular, one chooses a special  set of coset representatives
for $G$ modulo $G_0$, called Schreier representatives, and uses
these to find a set of generators for $G$. These generators are
labeled by the original generators of the group and the coset
representative.

\subsection{The relation between the action of the homeomorphism and
the surface kernel subgroup} \label{section:begin rewrite}

 We may assume that $S_0 = U/F_0$ where $F_0$ is the Fuchsian group
with presentation
\begin{equation} \label{equation:presentation}
\langle a_1, ..., a_{g_0}, b_1,..., b_{g_0}, x_1, ...x_t|\;
 x_1 \cdots x_t (\Pi_{i=1}^g[a_i,b_i]) =1;
 x_i^p=1\rangle.
\end{equation}

We summarize the result of \cite{G3, JGind} using facts about the
$n's$ and $m's$ defined in section \ref{section:mandn}. Let $\phi:
F_0 \rightarrow \mathbb{Z}_p $ be given by
$$\phi(a_i) = \phi(b_i) =0 \;\;\; \forall i=1,...,g_0 \mbox{ and }
\phi(x_j) = n_j \ne 0 \;\;\;\forall j=1,...,t.$$
If $F = Ker\; \phi$, then $S= U/F$. Moreover, $F_0/F$ acts on $S$
with quotient $S_0$. Conjugation by $x_1$ acts on $F$ and if $h$ is
the induced conformal map on $S$, $\langle h \rangle$ is isomorphic
to the action induced by this conjugation and the conjugacy class of
$h$ in the mapping-class group is determined by the set of $n_i$.
The order of the $n_i$'s does not affect the conjugacy class.
Replacing $h$ by a conjugate we may assume that $0 < n_i \le n_j <
p$ if $i < j$.

When we need to emphasize the relation of $h$ to $\phi$, we write
$h_{\phi}$ to mean the automorphism determined by conjugation by
$x_1$. The conjugacy class of $h^2$, would then be determined by the
homomorphism $\psi$ with $\psi(x_j) \equiv  2\phi(x_j) \; (p)$ or by
conjugation by $x_1^2$.

$F$ is sometimes called the surface kernel and $\phi$ the surface
kernel homomorphism \cite{GandG, Harvey}. Any other map from $F_0$
onto $\mathbb{Z}_p$ with the same $(m_1,....,m_{p-1})$ and with
$\phi(x_j) \ne 0 \;\; \forall j$ will yield an automorphism
conjugate to $h$.

\subsection{The rewriting} \label{section:rewrite}

We want to apply the rewriting process to words in the generators of
this  presentation for $F_0$ to obtain a presentation for $F$. We
choose right coset representatives for $H = \langle h \rangle$ as
$1, x_1,x_1^2,x_1^3,...,x_1^{p-1}$ and observe that this set of
elements form a {\sl Schreier} system. (see pages 88 and  93 of
\cite{MKS}). That is, that every initial segment of a representative
is again a representative.

 The Schreier right coset function assigns to a word $W$
in the generators of $F_0$, its coset representative $\overline{W}$
and $\overline{W} = x_1^q$ if $\phi(W)= \phi(x_1^q)$.

If ${\mathfrak{a}}$ is a generator of $F_0$, set
$S_{K,\;{\mathfrak{a}}} = K {\mathfrak{a}}
{\overline{K{\mathfrak{a}}}^{\;-1}}$. The rewriting process $\tau$
assigns to a word that is in the kernel of the map $\phi$, a word
written in the specific generators, $S_{K,\;{\mathfrak{a}}}$ for
$F$.
 Namely,  if ${\mathfrak{a}}_w$, $w = 1, ..., r$ are generators for
 $F_0$ and
$$U= {\mathfrak{a}}_{v_1}^{\epsilon_1}{\mathfrak{a}}_{v_2}^{\epsilon_2} \cdots
{\mathfrak{a}}_{v_r}^{\epsilon_r} \;\;\;\;\;(\epsilon_i = \pm 1),$$
defines an element of $F$, then (corollary 2.7.2 page 90 of
\cite{MKS})
$$\tau(U)= S_{K_1,\;{\mathfrak{a}_{v_1}}}
^{\epsilon_1} S_{K_2,\;{\mathfrak{a}}_{v_2}}^{\epsilon_2} \cdots
S_{K_r,\;{\mathfrak{a}}_{v_r}}^{\epsilon_r}$$ where $K_j$ is the
representative of the initial segment of $U$ preceding
${\mathfrak{a}}_{v_j}$ if $\epsilon_j = 1$ and $K_j$ is the coset
representative of $U$ up to and including
${\mathfrak{a}}_{v_j}^{-1}$ if $\epsilon_j = -1$.

In our case each ${\mathfrak{a}}_v$ stands for some generator of
$F_0$, that is one of the $a_i$ or $b_i$ or $x_j$.  We apply theorem
2.8 of \cite{MKS} to see
\begin{theorem} \label{theorem:Fpres} \cite{JGind} Let  $F_0$ have the presentation
given by equation \ref{equation:presentation}. Then $F$ has
presentation
\begin{eqnarray} \label{equation:presentationF}
\lefteqn{\langle S_{K,\;a_i}, \;S_{K,\;b_i},\; i = 1, ..., g_0;\;
S_{K,x_j}, \; j =
1,...t|}&\\
 && \tau(K \cdot x_1 \cdots x_t(\Pi_{i=1}^{g_0}[a_i,b_i]) \cdot  K^{-1})=1,\;  \tau(Kx_j^pK^{-1})=1 \rangle.
\end{eqnarray}
\end{theorem}

Here we  let $K$ run over a complete set of coset representative for
$\phi:F_0 \rightarrow H$ so that $F$ has generators
\begin{eqnarray*}
  S_{K,\;a_i}, S_{K,\;b_i} & &  i = 1, ..., g_0\\
  S_{K,\;x_j}, && j =1,...t
  \end{eqnarray*}

In \cite{JGind} we simplified the presentation and eliminated
generators and relations so that there is a single defining relation
for the subgroup. We first assumed that $\phi(x_1) = h$ noting that
if we can find a homology basis adapted to $h$, we can easily find a
homology basis adapted to any power of $h$ and, therefore, this
assumption will not be significant. We saw:

\begin{theorem}\label{theorem:homotopy} \cite{JGind} Let $F_0$ have the
presentation given by equation \ref{equation:presentation}. Then $F$
has \vskip .2in  $2pg_0 + (t-2)(p-1)$ generators: $$h^j(A_i),
h^j(B_i), i=1,....,g_0, j = 0, ..., p-1$$ $$ h^j(X_i), i=3,...,t, j
= 0, ..., p-2$$ and a single defining relation: $$  {\hat{\hat{R}}}
=1.  $$ Each generator and its inverse occur exactly once in
${\hat{\hat{R}}}$ and  every generator that appears is linked to
another distinct generator. ${\hat{\hat{R}}}$ is given explicitly in
terms of these generators.
\end{theorem}
\begin{corollary} \label{corollary:hom/hom} {\cite{JGind}} The homology basis obtained by abelianizing  the basis in theorem
\ref{theorem:homotopy} gives a homology basis adapted to $h$.
\end{corollary}

We will give ${\hat{\hat{R}}}$ explicitly in terms of the generators
once we have introduced more notation. In this paper we want to
simplify ${\hat{\hat{R}}}$. To do so we need to include some of the
calculations and notation from the earlier paper and illustrate
their use. We then introduce further notation and work to obtain the
desired result for homotopy (theorem \ref{theorem:BIG}).

\subsection{Illustration.}
We illustrate the use of  the notation:

 If we let $\phi(x_1) = h$ and $\phi(K)=\phi(x_1)^r$, then we have
$S_{K,X_j}= x_1^{r\phi(x_1)} \cdot x_j \cdot {\overline{K \cdot
x_j}}^{-1}$. Thus if
\begin{equation} \label{equation:proof}
 X_j =
x_j \cdot {\overline{x_j}}^{-1}, \mbox{ then } \;S_{K,x_j} =
h^r(X_j).
\end{equation}

We begin to rewrite the generators and relations using this
notation.

\noindent First we find $\tau(\overbrace{x_1x_1\cdots
x_1}^{p-\mbox{factors}})=1$. Since   $X_1= S_{1, x_1}$ we have
\begin{equation} \tau(x_1^p) = X_1 \cdot  h(X_1) \cdot h^2(X_1)
\cdots h^{p-2}(X_1)h^{p-1}(X_1)=1.
\end{equation}Similarly, if $\phi(x_j) = n_j$, and we set $X_j =
S_{{\overline{1}},x_j}$ and if ${\overline{K}} = x_1^{s}$, then we
can write $S_{K,\;x_j} = h^{s\cdot n_j}(X_j)$.

This tells us that
\begin{equation} \label{equation:homotopyXj}
\tau(x_j^p) = X_j \cdot h^{n_j}(X_j) \cdot h^{2n_j}(X_j) \cdots
h^{(p-2)n_j}(X_j) h^{(p-1)n_j}(X_j) =1 . \end{equation}

Note that in deriving all equations we are free to make use of the
fact (see \cite{MKS}) that $S_{M,\;x_1} \approx 1 \forall \mbox{
Schreier representatives } \; M$  where $\approx$ denotes {\sl
freely equal to}. Note that this also eliminates the $p$ generators, $S_{x_1^j, x_1}, j = 1,...,p$. 

Now  equation (\ref{equation:homotopyXj}) is a relation in the
fundamental group. We remind the reader that for a compact Riemann
surface, homology is abelianized homotopy so that when abelianized,
it reduces to
\begin{equation} \label{equation:homologyXj}
h^{p-1}(X_j) {\approx}^h -X_j - h(X_j) - \cdots - h^{p-2}(X_j)
\end{equation}  where
$\approx^h$ denotes {\it is homologous to}.

 Remember that our goal is to explicitly find the action on the
fundamental group given the conjugacy invariants in the mapping
class group. To this end we establish more notation to describe
${\hat{\hat{R}}}$ more simply and precisely.

\part{Homotopy and the algorithm}

\section{Finding the generators for the fundamental group and the algorithm} \label{section:LOOK}
While it is easy to see the action of $h$ on a homology basis,
writing explicitly the action of $h$ on generators for the
fundamental groups requires further notation. Our goal is to find a
{\sl homotopy basis} adapted to the automorphism $h$ and its
symplectic action on a symplectic basis obtained from this one.
\subsection{The action on the fundamental group of $S_0$}

First we note that when we rewrite $x_r^p=1$ as $\tau(x_r^p) = 1$,
if $\phi(x_r) = n_r$, then if $s_t\cdot  n_t \equiv n_1 \; (p)$ with
$0 \le s_t \le (p-1)$ and $\phi(x_1) = n_1$, then
\begin{equation} \label{equation:tauorder}
\tau(x_r^p) = S_{{\overline{1}},x_r}\cdot
S_{{\overline{x_r}},x_r}\cdot S_{{\overline{x_r^2}},x_r} \cdots
S_{{\overline{x_r^{p-2}}},x_r}  \cdot S_{{\overline{x_r^{p-1}}},x_r} = 1.
\end{equation}

For each $r$, this equation yields the homology relation of equation
\ref{equation:homologyXj} with $j$ replacing $r$.

If $\phi(x_j) \ne \phi(x_1)$ the relation (\ref{equation:tauorder})
(or equivalently (\ref{equation:homotopyXj})) determined by $x_j^p
=1 $ is slightly different than that obtained by replacing $j$ by
$1$. Up to homology both reduce to the same equation. However, the
statement regarding homotopy is more delicate.

 We note that if $\phi(x_1) = n_1$ and $\phi(x_r) =
n_r$ and $n_1 \cdot s_r \equiv n_r \;\;(p)$ and $s_r\cdot q_r \equiv
1 \;\; (p)$ so that $n_r \cdot q_r \equiv n_1 \;\; (p)$, then
$\phi(x_r)^{s_r} = \phi(x_1)$ and $\phi(x_r) = \phi(x_1^{q_r}) $ and
${\overline{x_r^k}} ={\overline{{x_1}^{k\cdot s_r}}} $.

If conjugation by $x_1$ induces the automorphism $h$ and $\phi(x_1)
=n_1 $ and $\phi(x_j) = n_j$ we let $n_1\cdot  s_j \equiv n_j
\;\;(p)$ so that conjugation by $x_j$ induces the same action on $F$
as conjugation by $x_1^{n_1\cdot s_j}$, and conjugation by $x_j^k$
induces the same action as conjugation by $x_1^{k\cdot  n_1\cdot
s_j}$. That is
\begin{equation} \label{equation:xjact}
S_{{\overline{x_j^k}},x_j}=
S_{{\overline{x_1^{k  \cdot s_j}}},x_j} =
h^{k \cdot s_j}(X_j)
\end{equation}

\subsection{Calculating with Simpler Notation}

In this section we first assume that $\phi(x_1)=1$ so that $h$ acts
as conjugation by $x_1$ and that $\phi(x_1) \le \phi(x_2) \le
\phi(x_j) \;\; \forall 3\le j \le t$. Then $s_j$  satisfies
$s_j\cdot n_j \equiv  1 \; (p)$.

To simplify the exposition from now on we assume that all exponents
are reduced modulo $p$ and are between $0$ and $p-1$. Instead of
working with the $X_j$ described above, we work with the $Y_j$
defined below. This is done so that the final result is in an
optimal form.

\vskip .1in
 We let $Y_j =
S_{\overline{x_1x_2 \cdots x_{j-1}},x_j}$ so that
$$h^z(Y_j) = S_{\overline{x_1^z\cdot x_1x_2 \cdots x_{j-1}},x_j},
\forall j=1,...,t \mbox{ and } \forall z = 0, \dots, p-1.$$ When
$K_r=\overline{ x_1x_2 \cdots x_{r-1}}$, (\ref{equation:tauorder})
gives \begin{align} \nonumber \tau(K_rx_r^pK_r) = & & \\
\nonumber S_{{\overline{x_1\cdots x_{r-1}}},x_r}\cdot
S_{{\overline{K_r\cdot x_r}},x_r}\cdot S_{{\overline{K_r\cdot
x_r^2}},x_r} \cdots &
  S_{{\overline{K_r\cdot x_r^{p-2}}},x_r}  \cdot S_{{\overline{K_r\cdot x_r^{p-1}}},x_r} = 1.
\end{align}
 or equivalently, where powers of $h$ are assumed to be reduced
modulo $p$,
\begin{equation} \label{equation:Ys}
Y_r \cdot h^{s_r}(Y_r) \cdot h^{2s_r}(Y_r) \cdot h^{(p-2)s_r}(Y_r) \cdot h^{(p-1)s_r}(Y_r) = 1.
\end{equation}

Note that the derivation of these  equations uses the fact that the
$S_{K,x_1} \approx 1$ for all coset representatives $K$.

 We will often drop the subscript
$r$ when it is understood and write (\ref{equation:Ys}) as
\begin{equation} \label{equation:nosubj}
Yh^s(Y)h^{2s}(Y) \cdots h^{(p-2)s}(Y)h^{(p-1)s}(Y) =1
\end{equation} where the powers of $h$ are again taken modulo $p$.

\begin{definition} {\it Notation:}  For any set of exponents
$u_1,...,u_v$ and any element $U$ of the group $F$,   we use the
shorthand notation $U^{u_1+u_2 + \cdots +u_v}$ and define
$$U^{u_1 + u_2 + \cdots +u_v} = h^{u_1}(U)h^{u_2}(U) \cdot h^{u_v}(U).$$
\end{definition}
Thus (\ref{equation:Ys}) becomes
\begin{equation} \label{equation:YYs} Y^{1 + s + 2s + \cdots (p-2)s + (p-1)s} =1.
\end{equation}
If we choose $p_r$ so that $1 \le p_r \le p-1$ and $p_r\cdot  s_r
\equiv p-1 \; (p)$, then we can solve (\ref{equation:Ys}) for
$h^{p-1}(Y_r)$ to obtain

\begin{equation}
\label{equation:hp-1Y} h^{p-1}(Y_r) =(Y_r^{(p_r+1)s_r + (p_r+2)s_r +
\cdots +(p-1)s_r} \cdot Y_r^{1 + s_r + 2s_r + \cdots
(p_r-1)s_r})^{-1}
\end{equation}

We now make some definitions. The motivation for these definitions
will become clear shortly.
\begin{definition} We define ${\bf{\mathcal{LR}}}_0$ by
$${{\bf {\mathcal{LR}}}}_0= (\Pi_{r=2}^t
Y_r^{p-1})(\Pi_{i=1}^{g_0}[h^{p-1}(A_i),h^{p-1}(B_i)])
$$
\end{definition}
 Substituting (\ref{equation:hp-1Y})
into this definition, we define
\begin{definition} We define
${\bf{\mathcal{LR}}}$ by

\begin{eqnarray*} \label{eqnarry:LR}
{{\bf {\mathcal{LR}}}} = & \\
& (\Pi_{r=2}^t
 (Y_r^{(p_r+1)s_r + (p_r+2)s_r + \cdots (p-1)s_r} \cdot
Y_r^{1 + s_r + 2s_r + \cdots+ (p_r-1)s_r})^{-1}) & \\
 &  \cdot (\Pi_{1=1}^{g_0}[h^{p-1}(A_i),h^{p-1}(B_i)])
& \\
\end{eqnarray*}
\end{definition}

\subsection{Eliminating generators and relations}
Recall that $h^r(Y_1) = S_{x_1^r,x_1} \approx 1$ for all integers
$r$. If $R$ is the relation in the group $$x_1\cdots
x_t\Pi_{i=1}^{g_0}[A_i,B_i]=1,$$ then with this notation we have:
$$\tau(R) = Y_2Y_3\cdots Y_t\Pi_{i=1}^g[A_i,B_i]=1  \mbox{  so that }
Y_2^{-1} = Y_3\cdots Y_t\Pi_{i=1}^g[A_i,B_i].$$

Let ${\bf \mathcal{Y}} = (Y_3\cdots Y_t\Pi_{i=1}^g[A_i,B_i])^{-1}.$

We have $\forall v = 1,..., p-1$
\begin{equation} \label{equation:Ysol}
\tau(x_1^vRx_1^{-v}) = h^v(Y_2)h^v(Y_3)\cdots
h^v(Y_t)\Pi_{i=1}^g[h^v(A_i),h^v(B_i)]=1 .
\end{equation}
Since $h$ is a group isomorphism, $h(TS) = h(S) h(T)$ so that
\begin{equation} \label{equation:caly}
h^v({\bf{\mathcal{Y}}})
 = (h^v(Y_3)\cdots h^v(Y_t)\Pi_{i=1}^g[h^v(A_i),h^v(B_i)])^{-1}.
 \end{equation}

Now $h^v(Y_2) = h^v({\mathcal{Y}})$ so we can eliminate   the
generator $Y_2$ and all of its images $h^v(Y_2)$ and replace these
by the expressions on the right hand side of equation
(\ref{equation:caly}). Next we replace the relation obtained from
equation (\ref{equation:Ys}) when $r=2$ with a relation in the
$h^v({\mathcal{Y}})$.

\begin{equation} \label{equation:nextYR}
{{\bf{\mathcal{Y}}}}^{1 + s_2 + 2s_2 + \cdots (p-1)s_2} =1.
\end{equation}

Recalling  that $p_2$ satisfies  $p_2 \times s_2 \equiv p-1 \; (p)$,
we have
\begin{equation} \label{equation:nextYRp2}
{{\bf{\mathcal{Y}}}}^{1 + s_2 + 2s_2 + \cdots +(p_2-1)s_2}
 \cdot {\bf{\mathcal{Y}}}^{(p-1)} \cdot
 {\bf{\mathcal{Y}}}^{(p_2 +1)s_2 +
(p_2+2)s_2 + \cdots + (p-1)s_2} =1.
\end{equation}

We observe that for each $j \ne 1,2$, every $h^{p-1}(Y_j)$ occurs in
$h^{p-1}({\bf{\mathcal{Y}}})$. Therefore, we eliminate these
$h^{p-1}(Y_j)$ using \ref{equation:hp-1Y}.

 Finally we recall the definition of
${\bf{\mathcal{LR}}}$, note that it is the same as
${\bf{\mathcal{Y}}}^{p-1}$ and substitute it into
(\ref{equation:nextYRp2}) to obtain

\begin{equation} \label{equation:finalR}
{{\bf{\mathcal{Y}}}}^{1 + s_2 + 2s_2 + \cdots +(p_r-2)s_2}
\cdot {\bf{\mathcal{LR}}} \cdot
 {\bf{\mathcal{Y}}}^{(p_r+1)s_2 +  \cdots +(p-1)s_2} =1
\end{equation}

This is now the single defining relation. We have proved

\begin{theorem}\label{theorem:BIG} {\bf Adapted Homotopy Basis.}
If the conjugacy class of $h$ in the mapping-class group is
determined by the $t$-tuple $(n_1,...,n_t)$ $t \ge 2$, $n_i \ne 0$,
$n_1=1 \le n_2 \le n_i, i\ge 3$, then $F$ has an {\sl adapted
homotopy basis}. That is, $F$ has  a presentation with \vskip .2in
\begin{equation} \label{equation:genssym}
{{\mbox{Generators:}}}
\;\;\;\;\;\;\;\;\;\;\;\;\;\;\;\;\;\;\;\;\;\;\;\;\;\;\;\;\;\;\;\;\;\;\;\;\;\;\;\;
\;\;\;\;\;\;\;\;\;\;\;\;\;\;\;\;\;\;\;\;\;\;\;\;\;
\end{equation}

\begin{eqnarray*} \label{eqnarry:genssym}
h^v(Y_w),&  w= 2, ..., t-2,&  v = 0,... ,p-2 \\
h^v(A_i), & i=1, \ldots, g_0, &v = 0,\ldots,p-1\\
h^v(B_i), &  i=1, \ldots, g_0, &v = 0, ..., p-1 \\
\end{eqnarray*}

and the single defining relation:

\begin{equation} \label{equation:surfacesymbol}
{{\bf{\mathcal{Y}}}}^{1 + s_2 + 2s_2 + \cdots +(p_2-2)s_2} \cdot
{\bf{\mathcal{LR}}} \cdot
 {\bf{\mathcal{Y}}}^{(p_2+1)s_2 +  \cdots +(p-1)s_2} =1,
 \end{equation}

where  ${\bf{\mathcal{Y}}}$ is the word in the generators:

$$(Y_3\cdots Y_t\cdot \Pi_{i=1}^{g_0}[A_i,B_i])^{-1},
$$

 ${\bf{\mathcal{LR}}}$ is the word in the generators:

\begin{align} \label{equation:word}
  (\Pi_{r=2}^t
 (Y_r^{(p_r+1)s_r + (p_r+2)s_r + \cdots (p-2)s_r}
Y_r^{1 + s_r + 2s_r + \cdots (p_r-1)s_r})^{-1}) \cdot \\
\nonumber  \Pi_{1=1}^{g_0}[h^{p-1}(A_i),h^{p-1}(B_i)],
\end{align} $s_r$ satisfies $s_r \cdot n_r \equiv 1\; (p)$ and
 $p_r$
satisfies $p_r \cdot s_r \equiv p-1 \; (p)$.

\end{theorem}

\begin{corollary} \label{corollary:symbol}
The fundamental group of $F$ has a surface symbol ${\mathcal{LR}}$
given by equation (\ref{equation:surfacesymbol}) whose unordered
edges are the generators given in equation (\ref{equation:genssym}).
\end{corollary}

We let ${\hat{\mathcal{LR}}}$ (respectively ${\hat{\mathcal{Y}}}$)
be the symbol ${\mathcal{LR}}$ (respectively ${\mathcal{Y}}$) from
which all of the occurrences of $h^v(A_i) \mbox{ and } h^v(B_i), i=
1, ...g_0, v = 0,...,(p-2)$ have been deleted.  We set $P =
\Pi_{i=1}^{g_0}[A_i,b_i]$ so that $h^v(P) =
\Pi_{1=1}^{g_0}[h^v(A_i),h^v(B_i)] $. Conjugating $h^v(A_i)$ and
$h^v(B_i)$  by the same word replaces $h^v(P)$ by an appropriate
conjugate. If we denote the appropriate conjugate by $[h^v(P)]_v$
and set ${\hat{\mathcal{Q}}} = P_0 \cdot  [h(P)]_1 \cdots
[h^{p-1}(P)]_{p-1}$ we can rewrite
\begin{equation} \label{equation:surfacesymbolnext}
{{\bf{\mathcal{Y}}}}^{1 + s_2 + 2s_2 + \cdots +(p_2-2)s_2} \cdot
{\bf{\mathcal{LR}}} \cdot
 {\bf{\mathcal{Y}}}^{(p_2+1)s_2 +  \cdots +(p-1)s_2} =1.
 \end{equation}
 as
\begin{equation}
{\hat{\mathcal{Q}}}  \cdot {{\bf{\hat{\mathcal{Y}}}}}^{1 + s_2 +
2s_2 + \cdots +(p_2-2)s_2} \cdot {\bf{\hat{\mathcal{LR}}}} \cdot
 {\bf{\hat{\mathcal{Y}}}}^{(p_2+1)s_2 +  \cdots +(p-1)s_2} =1,
\end{equation}
Setting  ${\hat{\hat{\mathcal{LR}}}} = {{\bf{\hat{\mathcal{Y}}}}}^{1
+ s_2 + 2s_2 + \cdots +(p_2-2)s_2} \cdot {\bf{\hat{\mathcal{LR}}}}
\cdot
 {\bf{\hat{\mathcal{Y}}}}^{(p_2+1)s_2 +  \cdots +(p-1)s_2}$ gives the surface symbol
${\hat{\mathcal{Q}}} \cdot {\hat{\hat{\mathcal{LR}}}}$.

With respect to this basis $h$ now has the matrix of the form
$$M_{{\tilde{\tilde{\mathcal{A}}}}} = \left(
\begin{array}{ccc}
M_{{\mathcal{A}}_{g_0, p\times p}} &  0 & 0\\
0 & M_{{\mathcal{A}}_{g_0, p\times p}} &  0 \\
0 & 0&  N_{{\hat{\hat{\mathcal{LR}}}}}
\\
\end{array} \right) $$
where $N_{{\hat{\hat{\mathcal{LR}}}}}$ is the matrix of $h$ with
respect to the ordered basis consisting of the edges of
${\hat{\hat{\mathcal{LR}}}}$.

\begin{corollary} \label{corollary:adsymb} The fundamental group has
the  evenly worded and fully and optimally linked surface symbol
${\hat{\mathcal{Q}}} \cdot {\hat{\hat{\mathcal{LR}}}}$ and matrix
$M_{\tilde{\tilde{{\mathcal{A}}}}}$.
\end{corollary}

\section{Applying the algorithm}

We now combine theorem \ref{theorem:BIG} and corollary
\ref{corollary:adsymb} with the tight linking algorithm of section
\ref{section:reducelink}. Initialize so that ${\mathcal{Q}} =
{\hat{\mathcal{Q}}}$, ${\mathcal{M}}=
M_{{\tilde{\tilde{\mathcal{A}}}}}$ and ${\mathcal{P}} =
{\hat{\hat{\mathcal{LR}}}}$ and conclude

\begin{theorem} \label{theorem:algbig}
Let $h$ be an element of the mapping-class group of prime order $p$
and assume that its  conjugacy class in the mapping-class group is
determined by the $t$-tuple $(n_1,...,n_t)$ $t \ge 2$, $n_i \ne 0$,
$n_1=1 \le n_2 \le n_i, i\ge 3$ so that  $F$ has the {\sl adapted
homotopy basis} described above. Then there is an algorithm that
inputs $(n_1,...,n_t)$ and outputs a unique symplectic matrix in the
conjugacy class of the image of $h$ in $Sp(2g,\mathbb{Z})$.
\end{theorem}

\begin{theorem} \label{theorem:algbignext}
Let $h$ be an element of the mapping-class group of prime order $p$
and assume that its  conjugacy class in the mapping-class group is
determined by the $t$-tuple $(n_1,...,n_t)$ $t \ge 2$, $n_i \ne 0$,
$n_1  \le n_2 \le n_i, i\ge 3$. Then there is an algorithm that
inputs $(n_1,...,n_t)$ and outputs a unique symplectic matrix in the
conjugacy class of the image of $h$ in $Sp(2g,\mathbb{Z})$.
\end{theorem}

\begin{proof} If $n_1 =1$, this is the previous theorem.
If $n_1 \ne 1$, let $s_1$ be the integer with $n_1 \times s_1 \equiv
1 \;\; (p)$. Consider the $t$-tuple whose $i$th entry is $n_i\times
s_1$ (reduced modulo $p$, of course). Apply theorem
(\ref{theorem:algbig}) to obtain a symplectic matrix and then output
this matrix raised to the power $n_1$.
\end{proof}

 We have also shown

\begin{theorem} \label{theorem:full}
 Consider an element $M$ of order $p$ in $Sp(2g,\mathbb{Z})$. For
each integer $t \in \{1,2,..., 2n+2\}$
 if there exists a set of $t$ integers
$\{n_1,...,n_t\}$ with $1 \le n_i \le (p-1)$ with $\Sigma_{i=1}^t
n_j \equiv 0 (p)$ or equivalently a $p-1$-tuple of real numbers
$(m_1,...,m_{p-1})$with $\Sigma_{i=1}^{p-1} im_i \equiv 0 (p)$ then
there is an element in the symplectic group of order $p$. The
integers $m_i$  determine its conjugacy class in the symplectic
group and the conjugacy class has a unique matrix representative
given by the algorithm. We call this matrix the normal form for $M$.
\end{theorem}
\begin{remark}
Note that there can be no $M$ of prime order for which $t=1$. The
canonical form for $t=0$ appears in section \ref{section:t=0}.
\end{remark}

Formulas for determining the number (possibly $0$) of distinct
$t$-tuples of integers that exist for any given $t$, $p$ and $g$
that satisfy the hypotheses of the theorem are to be found in
\cite{G3}.

\subsection{Complexity} \label{section:complex}
We note that replacing the complementary rotation numbers by the
rotation numbers requires an application of the Euclidean algorithm
$O(n^2)$. The calculation of the matrix
$M_{\tilde{\tilde{\mathcal{A}}}}$ is $O(g^2)$ because the words we
are looking at are of length at most $4g$. Computing the $pth$ power
of a $4g \times 4g$ matrix is $O(pg^3)$ but $p \le g$ so it is
$O(p^4)$. The complexity of the tight linking algorithm is $O(g^3)$.
Thus the complexity of the algorithm in theorem \ref{theorem:full}
is $O(g^4)$.

If $M$ is an $n \times n$ matrix, we define the norm of $M$ to be
the maximum of the absolute value of its entries and write $|M|$.
Multiplying matrices increases their norm and thus will affect the
running time. If we let $A_1,...,A_m$ be $m$,  $n \times n$
 matrices, then
$|A_1 \cdots A_m| =O(n^3 \cdot m \cdot \Sigma_{i=1}^m (1 + \log
|A_i|))$.

 We note in all applications
that the initial matrix ${\mathcal{M}}$ has entries that are only
$0$, $1$ and $-1$ and thus its norm is $1$. The matrix of the change
of basis at each iteration will also only have entries that are $0$,
$1$ and $-1$, so that the norm of it and its inverse is also $1$. We
can think of the final matrix as being obtained from the initial
matrix by conjugating it by at most $2g$ matrices so that $m=4g+1$.
Thus the norm of the final matrix in the tight linking algorithm is
$O(g^5)$. We may need to take up to $p$ powers of this matrix,
yielding a matrix of norm $O(g^{10})$.

\subsection{Other matrices} Let $M$ by a symplectic matrix of prime
order $p \ge 2$ with trace $T$. It follows from the results of
\cite{JGind}, that the map $\pi: MCG(S) \rightarrow
Sp(2g,{\mathbb{Z}})$ is injective on elements of finite order.
Either every preimage  of $M$ is an infinite order element of the
mapping-class group or there is a finite order preimage. In the
latter case, we have  $0 \le t \le 2g+2$ and the trace of
$M_{\tilde{A}} = 2-t$ so  since trace is a conjugacy invariant $T =
2-t$. Thus $M$ must be the image of an element of only infinite
order mapping-classes if $T < -2g$ or $T \ge 2$. There are prime
order elements of the symplectic group that are not the images of
finite order elements of the mapping-class group. If $-2g \le T \le
2$, then  a necessary condition for $M$ to have  a preimage of
finite order is that  there are integers $(n_1,...,n_t)$ such that
$\Sigma n_i \equiv 0 \;\; ({\mbox{mod }} p)$. Whether or not there
is such a t-tuple is determined in \cite{G3} where the number of
such $t-tuples$ is computed for each $t$ with $0 \le t \le 2g+2$.
The algorithm gives a sufficient condition, namely that $M$ be
conjugate to the unique normal form matrix.

\section{Homotopy when $t=0$} \label{section:t=0}

If the number of fixed points is zero, we can still find an adapted
basis. The calculations are slightly different. We have $2g =
2p(g_0-1) + 2$ and the presentation for the group $F_0$ is simply
\begin{equation} \label{equation:0presentation}
\langle a_1, ..., a_{g_0}, b_1,..., b_{g_0} |\;
  (\Pi_{i=1}^{g_0}[a_i,b_i]) =1;
 \rangle.
\end{equation}
Replacing $h$ by a conjugate if necessary, the map $\phi: F_0
\rightarrow \mathbb{Z}_p$ can be taken to be
$$\phi(a_i) = \phi(b_i) =0 \;\;\; \forall i=2,...,g_0 \mbox{ and }
\phi(a_1) = 1 \mbox{ and  } \phi(b_1) =0 $$ Using the rewriting with
coset representatives $1,a_1,...,a_1^{p-1}$, note that
$S_{a_1^{k},a_1} \approx 1$ for $k = 0, ... ,p-2$. Let $A =
S_{a_1^{p-1},a_1}$. Then $h$ acts on $Ker \; \phi$ via conjugation
by $a_1$.

We let $\{h^k(A_j),h^k(B_j), j = 2,...,g_0, k=0,...,p-1\}$ be as in
(\ref{equation:proof}) with $a_1$ playing the role of $x_1$ and let
$B= S_{1,b_1}$.  We let $P = \Pi_{i=2}^{g_0} [A_i,B_i]$.

Then we can compute that
\begin{align} \label{equation:new}
\tau(R) = 1 \implies h(B)B^{-1}P =1\\
\nonumber \tau(a_1^kRa_1^{-k}) = 1 \implies
h^k(B)(h^{k-1}(B))^{-1}h^k(P)
= 1 \forall k=1,...p-2\\
\nonumber  \mbox{and}\; \tau(a_1^{p-1}Ra_1^{-(p-1)})=1 \implies
ABA^{-1}(h^{p-1}(B))^{-1}h^{p-1}(P) = 1
\end{align}

We eliminate generators using (\ref{equation:new}) and let $\alpha =
A$ and $\beta = h^{p-1}(B)$ to obtain generators $$\{\alpha, \beta
\} \cup \{h^k(A_j),h^k(B_j), j = 2,...,g_0, k=0,...,p-1\}$$ and the
single defining relation $$\beta \alpha \beta^{-1}= h^{p-1}(P)
\alpha \Pi_{i=k}^{p-2}h^k(P).$$ Further we calculate that $h(\alpha)
= \alpha$ and $h(\beta) = \beta$. Note that $P$ is a product of
commutators. Thus the matrix representation for $h$ on $S= U/F$
where $F = \mbox{Ker} \phi$ is given by $2(g_0-1)$ permutation
matrices $M_{p \times p}$ and one two by two identity matrix. The
basis is a canonical homology basis. We have $\alpha \times \beta
=1$ and we note that the curves $\alpha$ and $\beta$ are by default
of type (1) in definition \ref{definition:adapt}

We have proved the  following:
\begin{theorem} {\label{theorem:0int}}
Assume $t=0$ and let $F$ be the  fundamental group of  the surface
$S$.
\begin{enumerate}
\item Then $F$
 \begin{enumerate}
  \item  has generators \begin{enumerate}
\item  $h^j(A_w), h^j(B_w)$  where $2 \le w \le g_0,  0 \le j \le p-1$.
\item $\alpha, \beta$
  where $h^k(\alpha) =  \alpha, h^k(\beta) = \beta, 0 \le k \le p$
\end{enumerate}
\item and the single defining relation
$$
 [\beta^{-1}, \alpha^{-1}] \cdot h^{p-1}(P^{\alpha}) \cdot
\Pi_{k=2}^{p-2}h^k(P) $$ where $P^{\alpha}$ denotes the words $P$
conjugated by $\alpha$ and  $P = \Pi_{i=2}^{g_0} [A_i,B_i]$.
\end{enumerate}
\item The corresponding intersection numbers are
\begin{enumerate}
\item $h^j(A_w) \times h^j(B_w) = 1$
\item $\alpha \times \beta =1$
\item  All other intersection numbers are $0$ except for those that
following from the above by applying the identities below to
arbitrary homology classes $C$ and $D$.
 \subitem $ C \times D = - D
\times C$ \subitem $ h^j(C) \times h^k(D)  = h^0(C) \times
h^{k-j}(D)$, ($k-j$ reduced modulo $p$.)
\end{enumerate}

\end{enumerate}
\end{theorem}

\section{Example, $p=3$, $t=5$, $(1,1,2,1,1)$} \label{section:ex}
In this section we elaborate on a specific example worked  out in
\cite{JGind}.  Assume $\phi(x_1) = 1 ,\phi(x_2) = 1,\phi(x_3) =
2,\phi(x_4) = 1,\phi(x_5) = 1$ so that $(n_1, ..., n_5) =
(1,1,2,1,1)$ and
 $(m_1,m_2)= (4,1)$.

 First replacing  $h$ by a conjugate,  we may assume that
 $\phi(x_1) = 1 ,\phi(x_2) = 1,\phi(x_3) =
1,\phi(x_4) = 1,\phi(x_5) = 2$. We choose as coset representatives
$x_1, x_1^2$ and $x_1^3$.

For any $g_0$, we have generators $$h^q(A_i),h^q(B_i), q= 0, \ldots,
p-1=2, i=1 \ldots, g_0.$$ and $$S_{{\overline{x_1^r}},x_j}, \;\;\; r
= 1,2,3, j = 2,3,4.$$ But $S_{{\overline{x_1^r}},x_1} \approx 1, r =
1,2,3$ and, therefore, these generators and the relation
$\tau(x_1^3)$ drops out of the set of generators and relations to be
considered.

We have
\begin{equation}
\tau(x_j^3) = S_{{\overline{1}}, x_j} \cdot S_{{\overline{x_j}},x_j}
\cdot S_{{\overline{x_j^2}},x_j}
\end{equation}

Set   $S_{{\overline{1}}, x_1}= Y_1$. Then $h(Y_1) = S_{x_1,x_1}$
and
 $h^2(X_j) = S_{x_1^2,x_j}$;

Set $S_{{\overline{x_1}}, x_2}= Y_2$ Then $h(Y_2) = S_{x_1^2,x_2}$
and
 $h^2(Y_2) = S_{x_1^3,x_2}$;

 Similarly, set
$S_{{\overline{x_1x_2}}, x_3}= Y_3$.  Then $h(Y_3) = S_{x_1^3,x_3}$
and
 $h^2(Y_3) = S_{x_1,x_3}$;
If $S_{{\overline{x_1x_2x_3}}, x_4}= Y_4$ Then $h(Y_4) =
S_{x_1,x_4}$ and
 $h^2(Y_4) = S_{x_1^2,x_4}$;and
finally if
 $S_{{\overline{x_1x_2x_3x_4}},
x_5}= Y_5$, then $h(Y_5) = S_{x_1^2,x_5}$ and
 $h^2(Y_5) = S_{x_1^3,x_j}$;

Use this notation and use $\tau(x_j^3)=1$ to see that

\begin{align} \label{equation:squares}
h^2(Y_2)\cdot Y_2 \cdot h(Y_2)= 1 \\
\nonumber h(Y_3)\cdot h^2(Y_3)\cdot Y_3=
1 \\
\nonumber Y_4\cdot h(Y_4)\cdot h^2(Y_4) = 1\\ \nonumber
h^2(Y_5)\cdot h(Y_5)\cdot Y_5 = 1 \end{align}

We compute \begin{equation} \label{equation:tauR} \tau(R) =
  S_{{\overline{1}},x_1}\cdot
 S_{{\overline{x_1}},x_2} \cdot
S_{{\overline{x_1x_2}},x_3} \cdot S_{{\overline{x_1x_2x_3},x_4}}
\cdot S_{{\overline{x_1x_2x_3x_4}},x_5} \cdot (\Pi_{i=1}^{g_0}
[A_i,B_i])=1.
\end{equation}
Using $S_{{\overline{x_1^r}},x_1} \approx 1, r = 1,2,3$ and solving
for $(Y_2)^{-1}$ in \ref{equation:tauR}, we have
\begin{equation}
(S_{x_1,x_2})^{-1} = Y_2^{-1} = Y_3 \cdot Y_4 \cdot Y_5  \cdot
(\Pi_{i=1}^{g_0} [A_i,B_i])=1.
\end{equation}
and
\begin{equation}
 (h(Y_2))^{-1} = h(Y_3) \cdot h(Y_4) \cdot h(Y_5)  \cdot
(\Pi_{i=1}^{g_0} [h(A_i),h(B_i)])=1.
\end{equation}
\begin{equation}
  (h^2(Y_2))^{-1} = h^2(Y_3) \cdot h^2(Y_4) \cdot h^2(Y_5)  \cdot
(\Pi_{i=1}^{g_0} [h^2(A_i),h^2(B_i)])=1.
\end{equation}
Using  $h^2(Y_2)Y_2h(Y_2) = 1$ and letting  $P=(\Pi_{i=1}^{g_0}
[A_i,B_i])$, we have
\begin{equation}
h(Y_3) \cdot h(Y_4) \cdot h(Y_5)  \cdot h(P)  \cdot Y_3 \cdot Y_4
\cdot Y_5  \cdot P \cdot h^2(Y_3) \cdot h^2(Y_4) \cdot h^2(Y_5)
\cdot h^2(P) =1
\end{equation}
We use equation  (\ref{equation:squares})  to replace the $h^2(Y_j)$
and
 obtain the relation

\begin{align} \label{equation:hathatR} h(Y_3) \cdot h(Y_4)  \cdot h(Y_5)  \cdot h(P)
  \cdot
Y_3 \cdot Y_4
\cdot Y_5  \cdot P \cdot (h(Y_3))^{-1} \cdot (Y_3)^{-1} & \\
 \nonumber \cdot (h(Y_4))^{-1} \cdot (Y_4)^{-1}   \cdot (Y_5)^{-1}
\cdot (h(Y_5))^{-1}\cdot h^2(P)\\ \nonumber &=1.
\end{align}

 Equation \ref{equation:hathatR} is the
relation ${\hat{\hat{R}}}=1$.

  Let $C_1 = h(Y_3) \cdot h(Y_4)
 \cdot h(Y_5)$,  $C_2 =  \cdot Y_3 \cdot Y_4
\cdot Y_5$ and \\
$C_3 =  (h(Y_3))^{-1} \cdot (Y_3)^{-1}(h(Y_4))^{-1} \cdot (Y_4)^{-1}
\cdot (Y_5)^{-1} \cdot (h(Y_5))^{-1}$

Let ${\tilde{P}} = h(P)^{{(C_2C_3)^{-1}}}\cdot P^{C_3^{-1}} \cdot
h^2(P)$ where if $X$ and $Y$ are words $X^Y$ means the conjugate of
$X$ by $Y$.

We can replace the relation \ref{equation:hathatR} by
\begin{eqnarray*} \label{eqnarray*:try}
  \tilde{P} h(Y_3) \cdot h(Y_4)
 \cdot h(Y_5)   \cdot Y_3 \cdot Y_4
\cdot Y_5  \cdot P \cdot (h(Y_3))^{-1} \cdot (Y_3)^{-1} \\
\cdot (h(Y_4))^{-1} \cdot (Y_4)^{-1}   \cdot (Y_5)^{-1} \cdot
(h(Y_5))^{-1}=1.
\end{eqnarray*}

We use the notation of section \ref{section:matrixforms} in
particular the definition of the matrix $B$ given at the end of that
section.  We can compute from the formulas for intersection numbers
in Theorem \ref{theorem:GPint} that the relevant part of the
intersection matrix, $B$, is  the $6 \times 6$ submatrix that gives
the intersection matrix for the curves in the basis given in the
order
$$X_{1,3}, h(X_{1,3}),X_{1,4},h(X_{1,4}),X_{2,1}, h(X_{2,1})$$ is
$$B= I_{\hat{\hat{R}}}= \left(
  \begin{array}{rrrrrr}
    0 & 1 & 1 & 0 & 1 &-1 \\
    -1 & 0 &-1 & 1 & 0 & 1 \\
    -1 & 1 & 0 & 1 & 1 & -1 \\
    0 & -1 & -1 & 0 & 0 & 1 \\
    -1 & 0 & -1 & 0& 0 & 0 \\
    1 & -1 & 1 & -1 & 0 & 0 \\
  \end{array}
\right).$$

That is, the matrix $J_{{\tilde{\mathcal{A}}}}$ breaks up into
blocks $\left(
  \begin{array}{ccc}
    0 & I_{pg_0} & 0 \\
    -I_{pg_0} & 0 & 0 \\
    0 & 0 & I_{\hat{\hat{R}}} \\
  \end{array}
\right)$

 We
now rearrange the relation. To simply the notation we let $a= Y_3$,
$b=Y_4$ and $c = Y_5$. So that the relation becomes
\begin{eqnarray*}
h(a)\cdot h(b) \cdot h(c) \cdot h(P) \cdot a \cdot b \cdot c \cdot P
\cdot (h(a))^{-1}a^{-1}(h(b))^{-1} \\
\cdot b^{-1} \cdot c^{-1} \cdot (h(c))^{-1}\cdot h^2(P) = 1.
\end{eqnarray*}

We can also make the simplifying assumption, replacing the elements
that occur in $P$, $h(P)$ and $h^2(P)$ by conjugates,  that we are
merely working with the symbol

\begin{equation}
h(a)\cdot h(b) \cdot h(c)  \cdot a \cdot b \cdot c  \cdot
(h(a))^{-1}a^{-1}(h(b))^{-1} \cdot b^{-1} \cdot c^{-1} \cdot
(h(c))^{-1} = 1. \end{equation}

We replace generators and relations using the algorithm as follows:

Let $M=h(a)\cdot W_1\cdot h(b) W_2\cdot (h(a))^{-1}$ where $W_1=
\lambda$, $W_2 = h(c) \cdot a \cdot b \cdot c$. Set $W_3 = a^{-1}$
and $W_4 = b^{-1} \cdot c^{-1} \cdot (h(c))^{-1}$.

Let $N= W_3W_2(h(a))^{-1}$. Then

\begin{align} \label{equation:rel}
  h(a)\cdot h(b) \cdot h(c) \cdot a \cdot b
\cdot c  \cdot (h(a))^{-1}a^{-1}(h(b))^{-1} \cdot b^{-1} \cdot
c^{-1}\cdot (h(c))^{-1} \\
\nonumber  = [M,N]W_3W_2W_1W_4 \;\;\;\;\;\;\;\;\;\; \;\;\;\;\;\;\;\;\;\;\;\;\;\;\;\\
\nonumber
 = [M,N] \cdot
a^{-1} \cdot h(c) a b c b^{-1} \cdot c^{-1} \cdot (h(c))^{-1}\\
\nonumber =1.
\end{align}

At this point one can proceed by inspection and let  $[b,c]^*$
denote the conjugate of $[b,c]$ by $a^{-1} \cdot h(c) \cdot a$ to
obtain

$$[M,N]\cdot [b,c]^* \cdot [a^{-1}, h(c)] =1. $$

However, to follow the algorithm carefully, we would set $\tilde{M}
= a^{-1}\cdot h(c) \cdot a$ and $\tilde{N}= b \cdot c \cdot b^{-1}
\cdot c^{-1} \cdot a$.

Then equation (\ref{equation:rel})  becomes $$[M,N]\cdot [\tilde{M},
\tilde{N}] \cdot [b,c] =1.$$

Thus the canonical homology basis is given by $$\{ h^j(A_i),
h^j(B_i) \}, i= 1,...,g_0, j = 0, ..., p-1 \cup \{M, N, \tilde{M},
\tilde{N}, b, c \}.$$

The reordered basis
 $\{M,
\tilde{M}, b, N, \tilde{N}, c \}$ has intersection matrix
$$\left(
                                               \begin{array}{rrrrrr}
                                                 0  & 0 & 0 & 1 & 0 & 0 \\
                                                 0 & 0 & 0 & 0 & 1 & 0 \\
                                                 0 & 0& 0& 0& 0& 1\\
                                                 -1 & 0& 0& 0& 0& 0\\
                                                 0 & -1 & 0& 0& 0& 0\\
                                                 0 & 0& -1& 0& 0& 0\\
                                                                                              \end{array}
                                             \right). $$

 We
can compute the action of $h$ on these last six elements of the
homology basis. First we note that $h(b) \approx^h M -  c - b -a -
h(c)$ and $h(a) \approx^h -N  + h(c)+ b + c$. Therefore,

$c \mapsto h(c) \approx^h \tilde{M}$

$h(c) \mapsto -c - h(c) \approx^h -c - \tilde{M}$

$a \mapsto h(a) \approx^h -N + h(c) + b +c  \approx^h -N + \tilde{M}
+ b + c$

$b \mapsto h(b) \approx^h M -c -b -a -h(c)\approx^h M -c - b
-{\tilde{N}} - \tilde{M}$

$M \mapsto h(M) \approx^h  -\tilde{M} - N$, and

$N \mapsto h(N) \approx^h -c + M -N $.

Thus the matrix of the action of $h$ with respect to the ordered
basis $M,\tilde{M},b,N, \tilde{N}, c$ is the submatrix we have been
seeking. Namely,

$$N_{symp{\tilde{\mathcal{A}}}} = \left(
  \begin{array}{rrrrrr}
    0&1&0&-1&0&0 \\
    0 & -1 & 0&1&0&-1 \\
    1&-1&-1&-0&-1&-1 \\
    1&0&0&-1&0& -1\\
    0&1&1&-1&0&-1 \\
    0&1&0&0&0&0 \\
 \end{array}
 \right)$$

One can verify that this $6 \times 6$ matrix really is a submatrix
of a symplectic matrix, as it should be.

\section{Acknowledgements} The author thanks Yair
Minsky and the Yale Mathematics Department for their hospitality and
support while some of this work was carried out. The author also
thanks Moira Chas for stimulating conversations and Rubi Rodriguez
for some helpful comments.

\end{document}